\newcommand{\Z}{\ensuremath{\mathbb{Z}}}           

\newcommand{\kar}{\mathrm{char}}                   

\renewcommand{\P}{\ensuremath{\mathbb{P}}}         
\renewcommand{\O}{\ensuremath{\mathcal{O}}}        
\renewcommand{\L}{\ensuremath{\mathcal{L}}}        
\newcommand{\I}{\ensuremath{\mathcal{I}}}          
\newcommand{\St}{\ensuremath{\mathrm{St}}}         
\newcommand{\Spec}{\mathrm{Spec}}                  
\newcommand{\Gr}{\mathrm{Gr}}                      
\newcommand{\Hom}{\ensuremath{\mathrm{Hom}}}       
\newcommand{\SHom}{\ensuremath{\mathcal{H}om}}     
\renewcommand{\H}{\ensuremath{\mathrm{H}}}         

\newcommand{\GL}{\mathrm{GL}}                      
\newcommand{\SL}{\mathrm{SL}}                      

\newcommand{\<}{\langle}                           
\renewcommand{\>}{\rangle}

\documentclass[a4paper]{amsart}

\usepackage{amssymb,euscript,amsmath}
\usepackage[all]{xy}

\CompileMatrices

\newtheorem{Theorem}{Theorem}
\newtheorem{Lemma}{Lemma}
\newtheorem{Proposition}{Proposition}
\newtheorem{Corollary}{Corollary}

\newtheorem{Remark}{Remark}

\newenvironment{Proof}{{\bf Proof:}}{$\square$}

\newcommand{\dX}{\mathrm{dX}}
\renewcommand{\u}{\mathfrak{u}}
\renewcommand{\b}{\mathfrak{b}}
\newcommand{\g}{\mathfrak{g}}
\newcommand{\p}{\mathfrak{p}}

\title[Frobenius splitting of cotangent bundles of flag varieties]{
Frobenius splitting of cotangent bundles of \\flag
varieties and geometry of nilpotent cones
}

\date{\today}

\author{Shrawan Kumar, Niels Lauritzen, Jesper Funch Thomsen}

\address{
Department of Mathematics\\University of North Carolina\\Chapel Hill, 
NC 27599-3250\\USA}
\email{kumar@math.unc.edu} 
\address{
Department of Mathematics\\University of North Carolina\\Chapel Hill, 
NC 27599-3250\\USA
}
\curraddr{
Matematisk Institut\\Aarhus Universitet\\Ny Munkegade\\ 
DK-8000 \AA rhus C\\Denmark
}
\email{niels@mi.aau.dk}

\address{
Matematisk Institut\\Aarhus Universitet\\Ny Munkegade\\ 
DK-8000 \AA rhus C\\Denmark
}
\email{funch@mi.aau.dk}

\thanks{The first author was supported in part by NSF grant no. DMS-9622887. 
The second author was supported by The Danish Research Council during the year 97/98 at
the University of North Carolina. An inspiring and hospitable environment is gratefully
acknowledged.
The third author thanks the University of North Carolina for its hospitality during his
visits. We thank H.~H.~Andersen, B.~Broer, J.~C.~Jantzen, B.~Kostant, O.~Mathieu, V.~Mehta, T.~R.~Ramadas, T.~Springer and 
W.~van der Kallen for valuable
discussions.
}

\dedicatory{To the memory of A.~Ramanathan}

\begin{document}

\maketitle

Let $G$ be a semisimple, simply connected algebraic group
over an algebraically closed field of prime characteristic $p>0$. Let
$U$ be the unipotent part of a Borel subgroup $B\subset G$ and
$\u$ the Lie algebra of $U$. Springer \cite{Springer} has 
shown for good primes, that there is a $B$-equivariant isomorphism
$U\rightarrow \u$, where $B$ acts through conjugation on $U$ and
through the adjoint action on $u$ (for $G=\SL_n$ one has the well known
equivariant
isomorphism $X\mapsto X-I$ between unipotent and nilpotent upper triangular 
matrices). Fix a good prime $p$.
Then there is an isomorphism of homogeneous bundles
$X=G\times^B U \rightarrow G\times^B \u$, where the latter can
be identified with the cotangent bundle $T^*(G/B)$ of $G/B$.

Motivated in part by 
\cite{LauritzenThomsen} we establish a link between the $G$-invariant form 
$\chi$ on the Steinberg module $\St=\H^0(G/B, (p-1)\rho)$ and Frobenius 
splittings \cite{MehtaRamanathan} of 
the cotangent bundle $T^*(G/B)$:
The representation
$\H^0(G/B, 2(p-1)\rho)$ is a quotient of the functions
$\H^0(X, \O_X)$ on $X$ (here
$\H^0(G/B, M)$ denotes the $G$-module induced from the 
$B$-module $M$ and $\rho$ half the sum of the roots $R^+$ opposite to the roots of $B$).
There is a natural map 
$$
\varphi:\St\otimes \St \rightarrow  \H^0(X, \O_X)
$$ 
such that the multiplication
$\mu:\St\otimes \St \rightarrow \H^0(G/B, 2(p-1)\rho)$
factors through the projection $\H^0(X, \O_X)\rightarrow \H^0(G/B,
2(p-1)\rho)$. Surprisingly the simple situation of \cite{LauritzenThomsen} generalizes
in that $\varphi(a\otimes b)$ is a Frobenius
splitting of $X$ if and only if $\chi(a\otimes b) = 1$ (if and only if $\mu(a\otimes b)$ is
a Frobenius splitting of $G/B$). 

Frobenius splitting of the cotangent bundle in this setup has a number
of nice consequences. By filtering differential forms via a morphism
to a suitable partial flag variety and using diagonality of the 
Hodge cohomology and Koszul resolutions, we obtain the vanishing theorem
$$
\H^i(G/B, S \u^*\otimes \lambda) = 0, i > 0
$$
where $\lambda$ is a dominant weight and $S\u^*$ denotes the symmetric
algebra of $\u^*$. This
was proved in \cite{AndersenJantzen} for large dominant weights and
for all dominant weights for groups of classical type and $G_2$ (and large primes).
The simple key lemma in the very simple proof
of the Borel-Bott-Weil theorem \cite{Demazure} implies that the vanishing theorem can be
extended to weights $\{\lambda\, |\, \<\lambda, \alpha^\vee\> \geq -1,\, \forall
\alpha\in R^+\}$. This vanishing 
theorem was proved 
in characteristic zero
by Broer \cite{Broer} using complete reducibility and the Borel-Bott-Weil theorem.
As in characteristic zero (\cite{Broer}, Theorem 4.4) it follows that the 
subregular nilpotent variety is normal, Gorenstein and has rational singularities.
In the parabolic case we prove the 
above vanishing theorem for regular dominant weights (after proving that the cotangent
bundle of partial flag varieties is also Frobenius split).

Another consequence is the conjectured isomorphism in (\cite{Jantzen}, 
II.12.15, \cite{AndersenJantzen}).
Furthermore by using the $B$-module
structure of $\St\otimes\St$, it follows easily that 
$T^*(G/B)$ carries a canonical Frobenius 
splitting \cite{Mathieu}\cite{Kallen}. This implies that 
$$
\H^0(G/B, S\u^*\otimes \lambda)
$$
has a good filtration \cite{Kallen} for any weight $\lambda$. One obtains (for all
groups in a uniform manner) that 
the 
cohomology of induced representations $\H^i(G_1, \H^0(G/B, \mu))^{[-1]}$ has a good filtration 
\cite{AndersenJantzen}(for
$\mu$ dominant and $p$ bigger than the Coxeter number of $G$).

Our canonical splitting relates to the splitting of Mehta and van der Kallen \cite{MehtaKallen} in
the $\GL_n$-case by taking a certain homogeneous component. For now we have ignored the more
combinatorial aspects of the methods in this paper, like analyzing compatible Frobenius
splitting.

\section{Notation and preliminaries}

The following notation is used throughout the paper.
Fix an algebraically closed field $k$ of characteristic $p>0$. All schemes and
morphisms will be over $k$. 

\subsection{Group data}
Let $G$ be a connected, simply connected semisimple algebraic group, $B$
a Borel subgroup of $G$, $T \subset B$ a maximal torus and $U$ the unipotent
radical of $B$. The Lie algebras of $G$, $B$ and $U$
are denoted $\g$, $\b$ and $\u$ respectively. 
In the following $B$ will act on $U$ by conjugation and on $\u$ by the
adjoint action. 
Let $B^+$ be the opposite Borel subgroup with
unipotent radical $U^+$, $R = R(T, G)$ the root system
of $G$ with respect to $T$, $R^- = R(T, U)$ (the negative roots), 
$R^+ = R(T, U^+) = \{\alpha_1, \dots, \alpha_N\}$ (the positive roots), $S\subset R^+$ the
simple roots and $h$ the Coxeter number
of $G$. 
For a parabolic subgroup $P\supset B$ we let $U_P$ denote the unipotent
radical of $P$, $U_P^+$ the opposite unipotent radical of $P$, $\mathfrak{u}_P$ 
the Lie algebra of $U_P$, $\p$ the Lie algebra of $P$ and $R_P\supset T$ the Levi factor of $P$.
By $\<\cdot, \cdot\>$ we denote the natural pairing $X(T)\times Y(T) \rightarrow \Z$ 
given by 
$\<\lambda, \mu\> = \lambda(\mu(1))$, 
where $X(T)$ is the group of characters (also identified with the weight lattice) and 
$Y(T)$ the group of one parameter subgroups of $T$ (also identified with the coroot lattice).
A simple root $\alpha\in R^+$ defines the (simple) reflection $s_{\alpha}(\lambda)=
\lambda - \<\lambda, \alpha^\vee\> \alpha$, where $\lambda\in X(T)$ and 
$\alpha^\vee\in Y(T)$ is the coroot associated with $\alpha$. 
For a subset $I\subset S$ we let $P=P_I$ denote the associated
parabolic subgroup. Recall that the group of characters $X(P)$ can be identified with 
$\{\lambda\in X(B) | \<\lambda, \alpha^\vee\> = 0, \text{for all\ }\alpha\in I\}$.
A weight $\lambda\in X(B)$ is called {\it dominant} if $\<\lambda, \alpha^\vee\> \geq 0$
for all $\alpha\in S$. A dominant weight $\lambda\in X(P)$ is called {\it $P$-regular} if 
$\<\lambda, \alpha^\vee\>>0$ for all $\alpha\not\in I$, where $P=P_I$ is a parabolic subgroup. 
A $B$-regular dominant weight is called {\it regular}.
The Weyl group $W$ of $G$ is generated by the simple
reflections.
The ``dot'' action of $W$ on $X(T)$
is given by $w\cdot \lambda = w(\lambda+\rho)-\rho$, where $\<\rho, \alpha^\vee\>=1$
for every simple root $\alpha\in S$.
On the weight lattice $X(T)$ the integral cone $\Z_+R^+\subseteq X(T)$ defines the
partial order: 
$ \lambda\geq \mu$ iff $\lambda-\mu\in \Z_+R^+$.

Recall that the prime $p$ is defined to be a {\it good prime} for $G$ if $p$
is coprime to all the coefficients of the highest root of $G$ written in terms
of the simple roots. For $G$ of almost simple type,
$p$ is a good prime if $p\geq 2$ for type $A$;
$p\geq 3$ for the types $B$, $C$ and $D$, $p\geq 5$ for the types
$F_4$, $E_6$, $E_7$ and $G_2$; $p\geq 7$ for the type $E_8$.

\subsection{Homogeneous bundles}

A $P$-scheme $X$ gives rise to an associated 
locally trivial fibration $G\times^P X$ over $G/P$ (\cite{Jantzen}, I.5.14, II.4.1). 
If $M$ is a finite dimensional $P$-representation, we let $\L(M)$ denote 
the sheaf of sections of the 
vector bundle $G\times^P M$ on $G/P$.

\subsection{The relative Frobenius morphism}

The {\it absolute Frobenius morphism} on a scheme is the identity on point
spaces and raising to the $p$-th power locally on functions. The absolute
Frobenius morphism is not a morphism of $k$-schemes.
Let 
$\pi:X\rightarrow \Spec(k)$ be a scheme.
Let $X'$ be the scheme obtained from $X$ by base change with the
absolute Frobenius morphism on $\Spec(k)$, i.e., the underlying topological
space of $X'$ is that of $X$ with the same structure sheaf $\O_X$ of rings, only
the underlying $k$-algebra structure on $\O_{X'}$ is twisted as 
$\lambda\odot f = \lambda^{1/p} f$, for $\lambda\in k$ and $f\in \O_{X'}$.
Using this description of $X'$, 
{\it the relative Frobenius
morphism} $F: X\rightarrow X'$ is defined in the same way as the 
absolute Frobenius morphism and it is a morphism of $k$-schemes.

\subsection{Frobenius splitting}

Following Mehta and Ramanathan \cite{MehtaRamanathan}
a variety $X$ is called
{\it Frobenius split} if the homomorphism $\O_{X'}\rightarrow F_*\O_X$ of
$\O_{X'}$-modules is split.
A homomorphism $\sigma:F_*\O_X\rightarrow \O_{X'}$ is a splitting
of $\O_{X'}\rightarrow F_*\O_X$
if and only if $\sigma(1)=1$. 
By abuse of terminology we will call an $\O_{X'}$-module homomorphism 
$\sigma:F_*\O_X\rightarrow \O_{X'}$ {\it a Frobenius splitting} 
if $\sigma(1)\in k\setminus \{0\}$ (so 
that $\sigma$ is a splitting up to a constant). 

A splitting $\sigma: F_*\O_X\rightarrow \O_{X'}$
is said to {\it split the subvariety $Y\subseteq X$ compatibly} if 
$\sigma(F_*\I_Y) \subseteq \I_{Y'}$, where $\I_Y$ denotes the ideal
sheaf of $Y$.

If $X$ is a smooth variety with canonical line bundle $\omega_X$, the Cartier operator gives
an isomorphism (\cite{MehtaRamanathan}, Proposition 5)
$$
\SHom_{\O_{X'}}(F_*\O_X, \O_{X'}) \cong F_* \omega_X^{1-p}.
$$
In this way global sections of $\omega_X^{1-p}$ correspond to homomorphisms
$F_*\O_X\rightarrow \O_{X'}$. A section of $\omega_X^{1-p}$ which corresponds
to a Frobenius splitting in this way, is called {\it a splitting section}. The above
isomorphism can be described quite explicitly in local coordinates (
\cite{MehtaRamanathan}, Proposition 5)

\begin{Proposition}
\label{PropositionLocal}
Let $P$ be a closed point of a smooth variety $Y$ over $k$ of dimension $n$. Choose
a system $x_1, \dots, x_n$ of regular parameters in the (regular) local ring 
$\O_{Y, P}$. Then the isomorphism
$$
F_*\omega_Y^{1-p} \rightarrow \SHom_{\O_{Y'}}(F_*\O_Y, \O_{Y'})
$$
is locally described as
$$
x^{\alpha}/(dx)^{p-1}: x^\beta\mapsto x^{((\alpha+\beta+\underline{1})/p)-\underline{1}},
$$
for any $\alpha=(\alpha_1, \dots, \alpha_n), \beta\in \Z_+^n$. Here we use the
multinomial notation $x^\alpha$ for the element $x_1^{\alpha_1}\dots x_n^{\alpha_n}
\in \O_{Y, P}$, and $\underline{m}=(m,\dots, m)\in \Z_+^n$ for an integer 
$m$. If $\gamma=(\gamma_1,\dots, \gamma_n)$ with at least one $\gamma_i$ nonintegral, 
we interpret $x^\gamma$ as zero. Furthermore $dx$ denotes the element
$dx_1\wedge \dots \wedge dx_n$, and $x^\alpha/(dx)^{p-1}$ denotes the local section
of $\omega_Y^{1-p}$ with value $x^\alpha$ on $(dx)^{p-1}$.
\end{Proposition}
\noindent
We also have
the following well known \cite{MehtaRamanathan}

\begin{Lemma}
\label{Lemmarestriction}
Let $U$ be an open dense subset of a smooth variety $X$. If a
section $s\in \H^0(X, \omega_X^{1-p})$ restricts to a splitting
section $s|_U\in \H^0(U, \omega_U^{1-p})$, then $s$ is a splitting
section.
\end{Lemma}

\begin{Lemma} 
\label{Lemmasplit}
Let $X$ be a Frobenius split variety and $\L$ a line bundle
on $X$. Then there is for each $i\geq 0$ an injection
$$
\H^i(X, \L) \hookrightarrow \H^i(X, \L^p)
$$
of abelian groups.
\end{Lemma}

\subsection{Volume forms}

Let $X$ be a smooth variety with trivial canonical bundle $\omega_X$. A {\it volume
form} is a nowhere vanishing section $\dX$ of $\omega_X$ (necessarily unique up to scalar
multiples if $\H^0(X, \O_X)^* = k$). {\it A function $f$ on $X$ is
said to Frobenius split $X$} (with respect to $\dX$) if 
$f\, \dX^{1-p}$ is a splitting section of $\omega_X^{1-p}$. 

\begin{Proposition}
\label{FsplitAn}
Let $X = \Spec\,k[x_1, \dots, x_n]$ be affine $n$-space. 
A volume form on $X$ is given
by $dx=dx_1\wedge \dots \wedge dx_n$ and a function $f\in k[X]$ Frobenius splits 
$X$ if and only if the coefficient of $x^{\underline{p-1}}$
in $f$ is nonzero and the coefficients of the terms $x^{\underline{p-1} + p\alpha}$
are zero for $\alpha\in \Z_{\geq 0}^n\setminus \{0\}$ (in the multinomial notation
of Proposition \ref{PropositionLocal}).
\end{Proposition}
\noindent
\begin{Proof}
An element $\sigma\in \Hom_{\O_{X'}}(F_*\O_X, \O_{X'})$ is a Frobenius splitting
if and only if $\sigma(1)$ is a nonzero constant. The proposition now follows
from Proposition \ref{PropositionLocal}.
\end{Proof}

\subsection{Filtration of differentials}

Let $f:X \rightarrow Y$ be a smooth morphism between smooth varieties $X$ and $Y$. Then
we have the following 

\begin{Lemma}
\label{Lemmafilt}
There is a short exact sequence
$$
0\rightarrow f^* \Omega_{Y/k} \rightarrow \Omega_{X/k} \rightarrow \Omega_{X/Y} \rightarrow 0,
$$
giving a natural filtration of $\Omega^m_{X/k}$ for $m\geq 1$ with associated graded object
$$
\Gr\, \Omega^m_{X/k} = \bigoplus_{i=0}^m f^* \Omega_{Y/k}^i\otimes \Omega_{X/Y}^{m-i}.
$$
\end{Lemma}

\subsection{The induction functor}
\label{SectionInduction}

Let $P$ be any parabolic subgroup.
For a $P$-module $M$ we let $\H^0(G/P, M)$ denote the induced $G$-module. Recall
that (in algebraic terms) $\H^0(G/P, M) = (M\otimes k[G])^P$, where $P$ acts
on $k[G]$ by right multiplication (it is a $G$-module with $G$ acting  trivially on $M$ and by
left multiplication on $k[G]$). This translates into the more familiar
$$
\H^0(G/P, M) = \{ f: G\rightarrow M | f(g\, p) = p^{-1}.f(g)\,\forall g\in G, p\in P \}.
$$
In this formulation $\H^0(G/P, M)$ is simply the global sections of the homogeneous
vector bundle $\L(M)$ on $G/P$. The sheaf cohomology 
$\H^i(G/P, \L(M))$ will also be denoted $\H^i(G/P, M)$ for $i\geq 0$.
For $P=B$, the functor $\H^0(G/B, -)$ is also 
denoted $\H^0(-)$.
If $M$ is a $G$-module, then $i:M\rightarrow \H^0(G/P, M)$ given by $i(m)(g) = g^{-1}.m$ is
an isomorphism of $G$-modules. 

\subsection{The Steinberg module}

The {\it Steinberg module} is by definition the induced module 
$\St = H^0(G/B, (p-1)\rho)$. It is irreducible
and selfdual. Fix an isomorphism $\St \rightarrow \St^*$ and denote the 
image of $v\in \St$ in $\St^*$ by $v^*$. This defines a $G$-invariant form given
by $\chi(v\otimes w) = \<v, w\> = v^*(w)$. Let $v^+$ and $v^-$ denote 
highest and lowest weight vectors of $\St$.

Let $G$ act on itself by conjugation. Then
the map $\St\otimes \St \rightarrow k[G]$ given by $(v\otimes w)(g) = \<v, g\, w\>$ 
is a $G$-homomorphism. We get in particular by restriction a $B$-homomorphism 
$$
\varphi: \St\otimes \St \rightarrow k[U].
$$
The global functions on $G\times^B U$ can be identified with
$\H^0(G/B, k[U])$. In this setting we have
$\H^0(\varphi)(v\otimes w)(g, u) = \<v, g u g^{-1} w\>$ using the identification
$i$ from \S \ref{SectionInduction}.

\subsection{The Frobenius kernel}

The relative Frobenius morphism $U\rightarrow U'$ is a homomorphism of
group schemes. The kernel $U_1$ is called the (first) {\it Frobenius kernel} and is a normal
(one point) subgroup scheme of $U$ (\cite{Jantzen}, I.9). If we fix a $T$-equivariant 
isomorphism (such that $x_i$ has weight $\alpha_i$)
$$
k[U] \rightarrow k[x_1, \dots, x_N],
$$ 
then $k[U_1] \cong k[x_1, \dots, x_N]/(x_1^p, \dots, x_N^p)$. 
Let $\gamma$ denote the
$B$-equivariant restriction homomorphism $k[U]\rightarrow k[U_1]$.
Notice 
that $k[U_1]$ is a finite dimensional $B$-representation with all
weights $\leq 2(p-1)\rho$. The $T$-equivariant projection on the highest weight space
spanned by $\bar{x}_1^{p-1}\dots \bar{x}_N^{p-1}$
is in fact a $B$-homomorphism $\psi: k[U_1]\rightarrow 2(p-1)\rho$, where the bar denotes the
corresponding element in $k[U_1]$.

\section{Frobenius splitting of $G\times^B U$}

We begin with the following elementary lemma.
\begin{Lemma}
\label{LemmaCanlinb}
For any parabolic subgroup $P$,
the canonical line bundle on the varieties $G\times^P U_P$ and $G\times^P \u_P$ is 
($G$-equivariantly) trivial.
\end{Lemma}
\begin{Proof}
We give the proof in the case $G\times^P U_P$. The argument for $G\times^P \u_P$ is
similar (in fact this is for good primes isomorphic to the cotangent bundle of
$G/P$).
Let $n=\dim U_P$.
The restriction of the locally free sheaf of relative differentials 
$\Omega = \Omega_{(G\times^P U_P)/(G/P)}$ on $G\times^P U_P$ to
$U_P = \{P\} \times^P U_P $ is the sheaf of differentials of $U_P$, and hence 
$\Omega^n|_{U_P} = \omega_{U_P}$. Let $dU_P$ be a volume form on $U_P$. Since $k[U_P]$
has no nonconstant units, the
canonical action of $P$ on $dU_P$ gives rise to a character $\beta$ of $P$, which
can be determined by considering the action of $P$ on $\omega_{U_P}|_e$, as the
identity $e\in U_P$ is fixed under $P$. The cotangent space at $e$ is
canonically isomorphic to $\mathfrak{M}_e/\mathfrak{M}_e^2$, where $\mathfrak{M}_e$
denotes the maximal ideal of functions in $k[U_P]$ vanishing at $e$. Hence
$\beta = \sum_{\alpha\in R(T, U^+_P)} \alpha$. Since $\Omega^n$ is a $G$-sheaf
it is the pull back of the line bundle induced by $\beta$ on $G/P$. As the
canonical line bundle on $G/P$ is induced by $-\beta$ the result follows from
Lemma \ref{Lemmafilt}.
\end{Proof}

Fix $T$-eigenfunctions $y_1, \dots, y_N\in k[U^+]$ of weights $-\alpha_1,\dots, -\alpha_N$, 
such that $k[U^+] \cong k[y_1, \dots, y_N]$. 
By Lemma \ref{LemmaCanlinb}, $X=G\times^B U$ carries a volume form $\dX$ 
restricting to $dy_1\wedge \dots \wedge dy_N \wedge dx_1 \wedge \dots \wedge dx_N$ on
the open subset $U^+ \times U \hookrightarrow G \times^B U$.
The following lemma is instrumental in proving Frobenius splitting of
$G\times^B U$.

\begin{Lemma}
\label{nonzero}
The map $\psi\circ \gamma\circ \varphi:\St\otimes \St\rightarrow 2(p-1)\rho$ is non-zero.
\end{Lemma}
\noindent
\begin{Proof}
We need to prove that the monomial $x_1^{p-1} \dots
x_N^{p-1}$ occurs with non-zero coefficient in $f\in k[U]$, where
$f(x) = \<v^+, x\, v^+\>$. The functions $x\mapsto \<v^+, x\, v^-\>$ and
$x\mapsto \<v^-, x\, v^-\>$ from $G$ to $k$ are highest and lowest weight vectors
in $\St = \H^0(G/B, (p-1)\rho)$ respectively. By Theorem 2.3 in \cite{LauritzenThomsen}
the function $\sigma$
$$
x\mapsto \<v^+, x\, v^-\> \<v^-, x\, v^-\> \in \H^0(G/B, 2(p-1)\rho)
$$
is a splitting section of $G/B$. The restriction of
$\sigma$ to $U^+$ is given by $x\mapsto \<v^-, x\, v^-\>$. Since $f$
corresponds to this function (which Frobenius splits $U^+$) under conjugation 
with $w_0$ (the longest element in $W$), the coefficient of 
$x_1^{p-1}\dots x_N^{p-1}$ in $f$ must be nonzero by Proposition \ref{FsplitAn}.
\end{Proof}

If $M$ is a $G$-module and $N$ a $B$-module, then by Frobenius reciprocity,
restriction followed by evaluation at $e\in G$ is an 
isomorphism (\cite{Jantzen}, Proposition I.3.4) 
$$
\Hom_G(M, \H^0(G/B, N)) \rightarrow \Hom_B(M, N).
$$
Let $\mu:\St\otimes \St \rightarrow \H^0(G/B, 2(p-1)\rho)$ denote the
multiplication map.

\begin{Corollary}
\label{comm}
There is a commutative diagram
$$
\xymatrix{
\H^0(G/B, k[U]) \ar[r]_{\H^0( \gamma)} & \H^0(G/B, k[U_1])\ar[d]^{\H^0( \psi)}\\
\St\otimes \St \ar[u]^{\H^0( \varphi)} \ar[r]^\mu & \H^0(G/B, 2(p-1)\rho)
}
$$
of $G$-equivariant homomorphisms.
\end{Corollary}
\noindent
\begin{Proof}
By applying the induction functor we get a homomorphism
$$
\H^0( \psi)\circ\H^0( \gamma)\circ \H^0( \varphi): \St\otimes \St
\rightarrow \H^0(G/B, 2(p-1)\rho)
$$
which is non-zero by Lemma \ref{nonzero} (and Frobenius reciprocity). By Frobenius reciprocity
$\mu$ is (up to a constant) the unique  
$G$-homomorphism $\mu:\St\otimes \St\rightarrow \H^0(G/B, 2(p-1)\rho)$). Adjusting constants this gives 
that the diagram is commutative. 
\end{Proof}

\begin{Theorem}
\label{Theoremmain}
Let $v=\sum_i v_i\otimes w_i$ be an element of $\St\otimes\St$.
The function $f=\H^0( \varphi)(v)$  
Frobenius splits $G\times^B U$ 
if and only if $\mu(v)$ is a splitting section of $\omega_{G/B}^{1-p}$.
In particular the function $f=f_v:G\times^B U\rightarrow k$
given by 
$$
f_v(g, u) = \sum_i \<v_i, g u g^{-1} w_i\>
$$ 
for $g\in G$, $u\in U$, Frobenius splits $G\times^B U$ if and 
only if $\chi(v)$ is nonzero. 
\end{Theorem}
\noindent
\begin{Proof}
Suppose that $\mu(v)$ is a splitting section of $\omega_{G/B}^{1-p}$. Let $f=\H^0( \varphi)(v)$. 
We prove that $f$ Frobenius splits $X=G\times^B U$ with respect to the volume form $dX$.
Restrict $f\,\dX^{1-p}$ to the open subset
$U^+\times U \hookrightarrow G\times^B U$. This leads to a 
form $f' (dy_1 \wedge \dots \wedge dy_N \wedge dx_1 \wedge \dots \wedge dx_N)^{1-p}$ on $U^+\times U$. 
By Proposition \ref{FsplitAn} and Lemma \ref{Lemmarestriction},
we are done if we prove that the monomial $y^{\underline{p-1}}
x^{\underline{p-1}}$ occurs with nonzero coefficient in $f'$ and 
the monomials $y^{\underline{p-1} + p \alpha} 
x^{\underline{p-1} + p \beta}$ occur with zero coefficient
where $\alpha, \beta\in \Z_{\geq 0}^N$ not simultaneously zero (in the multinomial
notation of Proposition \ref{PropositionLocal}).
We have the following commutative diagram 
$$
\xymatrix{
k[U]\otimes k[U^+] \ar[r]^{\gamma\otimes 1} & k[U_1] \otimes k[U^+]\ar[r]^{\psi\otimes 1} & 
2(p-1)\rho\otimes k[U^+]\\
(k[U]\otimes k[G])^B \ar[u]\ar[r]^{\H^0( \gamma)} 
& (k[U_1]\otimes k[G])^B \ar[u]\ar[r]^{\H^0( \psi)} & (2(p-1)\rho \otimes
k[G])^B \ar[u]\\
}
$$
with natural $T$-equivariant maps.
A monomial 
$y^{\underline{p-1} + p \alpha} 
x^{\underline{p-1} + p \beta}$ occuring in $f'$ must have
$\beta = 0$, as it is the restriction of an element in the image of
$(\St\otimes \St\otimes k[G])^B \rightarrow (k[U]\otimes k[G])^B$
and since a weight in $\St\otimes\St$ is $\leq 2(p-1)\rho$.
Furthermore by Corollary \ref{comm} 
$(\H^0( \psi)\circ \H^0( \gamma))(f)$ restricted to $U^+$ is a Frobenius splitting. Chasing 
through
the above diagram this means (using $\beta=0$) that $\alpha=0$ and the
monomial $y^{\underline{p-1}}
x^{\underline{p-1}}$ occurs with nonzero coefficient in $f'$, so that $f$  
Frobenius splits $G\times^B U$. On the other hand if $\H^0( \varphi)(m)$ is a 
Frobenius splitting it is easy to read off the diagram that $\mu(m)$ is a
splitting section. The last part of the theorem follows from Theorem 2.3 in 
\cite{LauritzenThomsen}.
\end{Proof}

Recall that the cotangent bundle $T^*(G/P)$ of $G/P$ is the $G$-fibration
associated to the $P$-module $(\g/\p)^*$ under
the adjoint action. It is well known that there is an
isomorphism $(\g/\p)^*\cong \u_P$ of $P$-modules in good
characteristics (\cite{Springer}, Lemma 4.4). Hence in this case 
$T^*(G/P)\cong G\times^P \u_P$.
We have the following crucial result due to Springer
(\cite{Springer}, Proposition 3.5)

\begin{Proposition}
\label{PropositionSpringer}
Let $\kar\, k$ be a good prime for $G$. Then there exists a $B$-equivariant
isomorphism $\zeta: U\rightarrow \mathfrak{u}$. Moreover for any parabolic subgroup
$P$, $\zeta$ restricts to give a $P$-equivariant isomorphism $\zeta_P: U_P\rightarrow
\mathfrak{u}_P$.
\end{Proposition}

\begin{Corollary}
\label{Corollaryfsplit}
Let $\kar\, k$ be a good prime for $G$. Then
the cotangent bundle $T^*(G/B)$ of $G/B$ is Frobenius split.
\end{Corollary}
\noindent
\begin{Proof}
By Proposition \ref{PropositionSpringer} we get a $G$-isomorphism 
$G\times^B U \rightarrow G\times^B \u$, where the latter can be identified 
with the cotangent bundle of $G/B$. The result
now follows from Theorem \ref{Theoremmain}.
\end{Proof} 

\begin{Remark}
For $v\in \St\otimes \St$ define $f_v:G\times^B B\rightarrow k$ as in Theorem 
\ref{Theoremmain}. Then $f_v$ Frobenius splits $G\times^B B$ if and
only if $\chi(v)\neq 0$. Furthermore any such $v$ gives rise to a Frobenius
splitting of $G\times^B \b$, which descends via the map $(g, X)\mapsto
Ad(g)X$ to the Lie algebra $\g$. Since we have no nontrivial applications of
these results we do not give any proofs.
\end{Remark}

\section{Vanishing}
\noindent
Let 
$$
\mathcal{C}=\{\mu\in X(T) | \<\mu, \alpha^\vee\> \geq -1, \forall \alpha\in R^+\}.
$$
It is easy to see (\cite{Broer1}, Proposition 2) that $\mathcal{C}$ is the set of 
weights $\lambda$ such that
if $\mu$ is a dominant weight with $\lambda \leq \mu \leq \lambda^+$, 
then $\mu = \lambda^+$ (here $\lambda^+$ denotes the dominant weight in the
$W$-orbit of $\lambda$). The set $\mathcal{C}$ is precisely the weights of line bundles
on $G/B$ in characteristic zero, 
which have vanishing higher cohomology when pulled back to
the cotangent bundle (\cite{Broer}, Theorem 2.4).
In this section we prove the analogous vanishing theorem in good prime
characteristics. Andersen and Jantzen (\cite{AndersenJantzen}, Theorem 3.6)
proved the following vanishing theorem under the assumption that $\lambda=0$ or $\lambda$
strongly dominant ($\<\lambda, \alpha^\vee\>\geq h-1$ for all $\alpha\in S$).
For $p\geq h-1$ and all components of $G$ classical or $G_2$ 
they proved the vanishing theorem for $\lambda$ dominant (\cite{AndersenJantzen}, 
Proposition 5.4). In fact the condition 
$\lambda+\rho$ dominant in (\cite{AndersenJantzen}, Proposition 5.4) is not sufficient 
for vanishing as noticed by Graham and Broer - this is also revealed using 
Lemma \ref{Lemmasimple} in \S \ref{Sectionsimplelemma} coupled with Bott's theorem.
Let $\pi: T^*(G/B)\rightarrow G/B$ denote the projection.
\begin{Theorem}
\label{Theoremvanishing}
Let $\kar\, k$ be a good prime for $G$ and
suppose that $\lambda\in \mathcal{C}$. Then
$$
\H^i(T^*(G/B), \pi^*\L(\lambda)) = \H^i(G/B, S u^*\otimes \lambda) = 0
$$
when $i>0$.
\end{Theorem}

\begin{Remark}
By the semicontinuity theorem our result implies the same vanishing theorem
over fields of characteristic zero.
\end{Remark}

\subsection{The Koszul resolution}
\label{SectionKoszul}
Let 
$$
0\rightarrow V' \rightarrow V \rightarrow V'' \rightarrow 0
$$
be a short exact sequence of vector spaces. For any $n>0$ one obtains a functorial
exact sequence (called the {\it Koszul resolution})
$$
\dots\rightarrow S^{n-i} V\otimes \wedge^i V'\rightarrow \dots \rightarrow
S^{n-1} V\otimes V'\rightarrow S^n V \rightarrow S^n V'' \rightarrow 0.
$$

\subsection{A simple lemma}
\label{Sectionsimplelemma}
Let $P_\alpha$ be the minimal parabolic subgroup corresponding to the
simple root $\alpha$. If $\lambda\in X(T)$ is a weight with $\<\lambda, 
\alpha\> = -1$ and $V$ a $P_\alpha$-module, then
$$
\H^i(G/B, V\otimes \lambda) = 0
$$
for $i\geq 0$. This result is the simple key lemma in Demazures very
simple proof of the Borel-Bott-Weil theorem \cite{Demazure}. It has the following
consequence (a similar approach has been used by Broer in \cite{Broer2})

\begin{Lemma}
\label{Lemmasimple}
Suppose that $\lambda\in \mathcal{C}$ and $\<\lambda, \alpha^\vee\> = -1$ for
a simple root $\alpha$. Then $s_\alpha(\lambda)\in \mathcal{C}$ and
$$
\H^i(G/B, S^n {\mathfrak u}^* \otimes \lambda) \cong 
\H^i(G/B, S^{n-1}{\mathfrak u}^* \otimes s_{\alpha}(\lambda))
$$
for $i\geq 0$ and $n>0$.
\end{Lemma}
\noindent
\begin{Proof}
As $s_\alpha$ permutes $R^+\setminus\{\alpha\}$ and maps $\alpha$ to
$-\alpha$, we get that $s_\alpha(\lambda)\in \mathcal{C}$.
The isomorphism follows by applying \S\ref{SectionKoszul} to the short
exact sequence of $B$-modules
$$
0\rightarrow \alpha \rightarrow \u^*
\rightarrow \u^*_{P_\alpha}\rightarrow 0,
$$
and then tensoring with $\lambda$,
where $P_\alpha$ is the minimal parabolic subgroup corresponding to the simple
root $\alpha$.
\end{Proof}

\subsection{Large dominant weights}

This section contains a proof of a lemma enabling us to turn
Frobenius splitting into vanishing for weights, which are
not necessarily regular. The key lies in filtering differentials using
the fibration $G/B \rightarrow G/P$ for a suitable parabolic subgroup 
$P\supset B$.

\begin{Lemma}
\label{Lemmalarge}
Let $\lambda$ be a dominant weight. Then
$$
\H^i(G/B, \Omega_{G/B}^j\otimes \L(m\lambda)) = 0
$$
for $i>j$ and $m$ sufficiently big.
\end{Lemma}
\noindent
\begin{Proof}
If $\lambda = 0$, 
we are done by the fact that $H^i(G/B, \Omega_{G/B}^j) = 0$ 
for $i\neq j$ (\cite{Jantzen}, II.6.18). This is usually referred to
as diagonality of Hodge cohomology.
If $\lambda\neq 0$, we can
choose a (unique) parabolic subgroup $P\neq G$, such that $\lambda$ is a ($P$-regular) 
character of $P$ and the induced line
bundle $\L(\lambda)$ is ample on $G/P$. Let $f$ denote the smooth $(P/B)$-fibration
$G/B\rightarrow G/P$. Using Lemma \ref{Lemmafilt}, we see that it is
enough to prove that the cohomology groups
$$
\H^i(G/B, f^*\Omega_{G/P}^r\otimes \Omega_{(G/B)/(G/P)}^{j-r}\otimes \L(m\lambda))
$$
vanish for all sufficiently big $m$, where $0\leq r\leq j$.
The $E_2$-terms in the Leray spectral sequence for $f$ are (using the projection 
formula)
\begin{align*}
E_2^{pq} &= \H^p(G/P, \L(m\lambda)\otimes \Omega_{G/P}^r\otimes R^q f_*\Omega_{G/B/G/P}^{j-r}) \\
         &= \H^p(G/P, \L(m\lambda)\otimes \Omega_{G/P}^r\otimes \L(H^q(P/B, \Omega_{P/B}^{j-r}))).  
\end{align*}
For all $m$ sufficiently big we get $E_2^{pq}=0$ for $p>0$ by Serre vanishing.
Diagonality of Hodge cohomology for $P/B$ gives that $E_2^{pq}=0$ unless $q=j-r$. 
In particular, for $m$ sufficiently big, combining the two, we get $E_2^{pq}=0$
unless $p=0$ and $q=j-r$. Now the result follows by the Leray spectral sequence, since
$i>j$ by assumption.
\end{Proof}

\subsection{Proof of theorem \ref{Theoremvanishing}}
\label{Proofoftheorem}
The first isomorphism follows since $\pi: T^*(G/B) \rightarrow G/B$ is an affine
morphism and $\pi_*\O_{T^*(G/B)} = \L(S \mathfrak{u}^*)$. To prove the
vanishing theorem we may assume that $\lambda$ is dominant, because of the
following argument: Assume by
induction on $n$ that $\H^i(G/B, S^j\mathfrak{u}^*\otimes\lambda)=0$ for $j<n$,
$i>0$ and $\lambda\in \mathcal{C}$. We wish to prove the same result 
for $j=n$. Take a non dominant weight $\lambda\in\mathcal{C}$.
Then there is a simple root $\alpha$ such that
$\<\lambda, \alpha^\vee\>=-1$. By Lemma \ref{Lemmasimple}, $s_\alpha(\lambda)\in
\mathcal{C}$ and
$$
\H^i(G/B, S^n\mathfrak{u}^*\otimes\lambda) =
\H^i(G/B, S^{n-1}\mathfrak{u}^*\otimes s_{\alpha}(\lambda)),
$$
where the latter group vanishes by induction. 

So assume that $\lambda$ is dominant.
Since $({\mathfrak b} / {\mathfrak u})^*$
is a trivial $B$-module, it follows from \S\ref{SectionKoszul} (applied to the
sequence $0\rightarrow (\mathfrak{b}/\mathfrak{u})^*\rightarrow \mathfrak{b}^* 
\rightarrow \mathfrak{u}^*\rightarrow 0$, and breaking the resulting Koszul
resolution up into short exact sequences) that the vanishing of
$\H^i(G/B, S \b^*\otimes \lambda)$ implies the vanishing of $\H^i(G/B, S \u^*\otimes \lambda)
$
for $i>0$.
Again using \S\ref{SectionKoszul} for the short exact sequence
$
0\rightarrow (\mathfrak g / \mathfrak b)^* \rightarrow 
{\mathfrak g}^* \rightarrow 
{\mathfrak b}^* \rightarrow 0
$
we get for $n\geq 1$ an exact sequence
$$
\dots
\rightarrow \wedge^1(\mathfrak{g}/\mathfrak{b})^*\otimes S^{n-1}\mathfrak{g}^*
\otimes\lambda \rightarrow S^n\mathfrak{g}^*\otimes\lambda\rightarrow
S^n\mathfrak{b}^*\otimes\lambda \rightarrow 0
$$
after tensoring with $\lambda$. By breaking this up into short exact sequences, we see that 
the vanishing $\H^i(G/B, S \mathfrak{b}^*\otimes \lambda)=0$ for any fixed $i>0$ 
follows from the vanishing 
$$
\H^{i+j}(G/B, \wedge^j(\mathfrak g/\mathfrak b)^* \otimes \lambda) = 0
$$
for all $j\geq 0$. The $B$-representation 
$\wedge^j(\mathfrak g/\mathfrak b)^*$
induces the bundle of $j$-forms $\Omega_{G/B}^j$ on $G/B$. By Lemma \ref{Lemmalarge},
we get for all large enough $r$ that $H^{i+j}(G/B, \wedge^j(\mathfrak{g}/\mathfrak{b})^*
\otimes(p^r \lambda)) = 0$ for $j\geq 0$ and hence $\H^i(G/B, S\mathfrak{u}^*\otimes 
(p^r \lambda))=0$ for $i>0$. But
by Corollary \ref{Corollaryfsplit} and
Lemma \ref{Lemmasplit}, we have an injection of abelian groups
$$
\H^i(T^*(G/B), \pi^*\L(\lambda))\hookrightarrow 
\H^i(T^*(G/B), \pi^*\L(p^r\lambda)) 
$$
which translates into an injection $\H^i(G/B, S \mathfrak{u}^*\otimes\lambda)
\hookrightarrow \H^i(G/B, S \mathfrak{u}^* \otimes (p^r\lambda))$
for any $r>0$ (this is where the assumption that $p$ is good for $G$ is used). This proves 
the theorem.

\subsection{Dolbeault vanishing}

Theorem \ref{Theoremvanishing} is in fact equivalent to the following (Dolbeault) vanishing
(see \cite{Broer1} for results in characteristic zero and the parabolic case).

\begin{Theorem}
Let $\kar\, k$ be a good prime for $G$ and
$\lambda\in \mathcal{C}$. Then
$$
\H^i(G/B, \Omega_{G/B}^j\otimes\L(\lambda)) = 0
$$ 
for $i>j$.
\end{Theorem}
\noindent
\begin{Proof}
Theorem \ref{Theoremvanishing} implies that $\H^i(G/B, S^n\b^*\otimes
\lambda)=0$ for $i>0$, using induction on $n$ in the Koszul resolution (tensored with $\lambda$) 
coming from the short exact sequence 
$0\rightarrow (\b/\u)^*\rightarrow \b^*\rightarrow \u^*\rightarrow 0$. This vanishing now
fits in a similar induction on $n$ in the Koszul resolution (tensored with $\lambda$) coming from
the short sequence sequence $0\rightarrow (\g/\b)^*\rightarrow \g^*\rightarrow \b^*\rightarrow 0$. This 
gives the desired vanishing.
\end{Proof}

\section{The parabolic case}

In this section we prove that the cotangent bundle of $G/P$, where $P$ is
a parabolic subgroup is Frobenius split when $\kar\, k$ is a good prime
for $G$.

\subsection{Frobenius splitting of $G\times^P U_P$}
 
Let $P = P_I \supset B$ be a parabolic subgroup given by the subset $I\subset S$. 
Let $R_I$ denote the root system generated by $I$. The functions $k[(U_P)_1]$ on the Frobenius kernel of
$U_P$ is a finite dimensional $P$-representation of highest weight $(p-1)\delta_P$, where
$\delta_P = \sum_{\alpha\in R^+\setminus R_I^+}\alpha \in X(P)$. Observe that $-\delta_P$ is 
the weight inducing the canonical line bundle of $G/P$.
The canonical line bundle on $G\times^P U_P$ is
trivial by Lemma \ref{LemmaCanlinb}. The global functions $k[G\times^P U_P]$ can
be identified with $\H^0(G/P, k[U_P])$, which is naturally isomorphic to
$\H^0(G/B, k[U_P])=k[G\times^B U_P]$.
As in the case of a Borel subgroup $B\subset P$ we have a
natural $P$-equivariant map $\varphi_P:\St\otimes\St \rightarrow k[U_P]$.
The natural map
$$
\H^0(G/B, k[U_P])\rightarrow \H^0(G/B, k[(U_P)_1])\rightarrow
\H^0(G/P, (p-1)\delta_P),
$$
composed with $\H^0(G/P, \varphi_P)$ gives a map $\mu_P:\St\otimes\St\rightarrow 
\H^0(G/P, (p-1)\delta_P)$.

\begin{Theorem}
\label{Theoremparabolic}
Let $v=\sum_i v_i\otimes w_i$ be an element of $\St\otimes\St$.
The function $f=\H^0(G/P, \varphi_P)(v)$  
Frobenius splits $G\times^P U_P$ 
if and only if $\mu_P(v)$ is a splitting section of $\omega_{G/P}^{1-p}$.
The function $f=f_v:G\times^P U_P\rightarrow k$
given by 
$$
f_v(g, u) = \sum_i \<v_i, g u g^{-1} w_i\>
$$ 
for $g\in G$, $u\in U_P$, Frobenius splits $G\times^P U_P$ if and 
only if $\chi(v)$ is nonzero. 
\end{Theorem}
\noindent
\begin{Proof}
It follows by analogous
weight considerations for the restriction to $U_P^+\times U_P$ as
in the $B$-case, that $v\in \St\otimes\St$ maps
to a Frobenius splitting function of $G\times^P U_P$ if and only if
$\mu_P(v)$ is a splitting section of $G/P$ (a useful fact
here is that $\alpha\in R^+\setminus R_I^+$ contains a simple root outside $I$
with nonzero coefficient when written as a sum of simple roots).
In order to prove the last part of the theorem, we need to exhibit an element
$w\in \St\otimes\St$ such that $\H^0(\varphi_P)(w)$ is a Frobenius splitting (because
this implies that $\mu_P(w)$ is a Frobenius splitting, so that $\mu_P$ followed
by the evaluation map \cite{LauritzenThomsen} $\H^0(G/P, (p-1)\delta_P)\rightarrow k$
is a non-zero $G$-homommorphism $\St\otimes \St\rightarrow k$).

As proved in Theorem \ref{Theoremmain}, the function defined by 
$f(g,u) = \<v^-, g u g^{-1} v^+ \>$,
Frobenius splits $G\times^B U$. The restriction of this
function to $U^+ \times U$ therefore Frobenius splits $U^+ \times U$.
Observe that this restriction is given by 
$$ f(g,u) = \<v^-, g u v^+\>, \, g \in U^+, u \in U.$$
Let $w'_0$ be the longest element of the Weyl group of $R_P$ and let 
$v_0^+ = w'_0 v^+, v_0^- = w'_0 v^-$. Index the set of positive roots $\{\alpha_1, 
\dots, \alpha_N\}$ 
in such a manner that the first $n$ roots are the positive roots of $R_P$. Let $y_i:k \to U $ 
(resp.  $x_i:k \to U^+ $) be the root homomorphism corresponding to the 
root $-\alpha_i$ (resp. $\alpha_i$). 

Write $u = y_{N}(t_N) \dots y_{1}(t_1)$ 
and $g=x_{1}(s_1) \dots X_{N}(s_N)$.
Then 
$$u v^+ = y_{N}(t_N) \dots y_{n+1}(t_{n+1})
( \sum_{{l} \neq \underbar{p-1}}\,c_l t_n^{l_n} \dots t_1^{l_l}
v_{l}+ c t_n^{p-1} \dots t_1^{p-1} v^+_0),$$
where $v_{l}$ are weight vectors in $\St$ and $l=(l_1, \dots, l_n)$.
As $f$  Frobenius splits  $U^+ \times U$, we see that  
the coefficient of $t_{n+1}^{p-1} \dots 
t_N^{p-1}s_N^{p-1} \dots s_{n+1}^
{p-1}$ in 
$$\<v_0^-, x_{n+1}(s_{n+1}) \dots x_{N}(s_{N}) 
y_{N}(t_N) \dots y_{n+1}(t_{n+1}) v^+_0\> $$
is nonzero. By weight considerations, it therefore easily follows
that the function 
$$ f' : U_P^+ \times U_P \rightarrow k ,\,\,\, f'(g,u) = \<v_0^-, 
g u v^+_0\>$$
Frobenius splits $U_P^+ \times U_P$. But $f'$ extends to the function 
(again denoted by) $f' : G\times^P U_P \to k $ given by  
 $(g,u) \mapsto \<v_0^-, g u g^{-1} v^+_0\>$. We claim that this
Frobenius splits $G\times^P U_P$.  To see this, it suffices to observe 
that $U_P^+$ fixes $v_0^+$.
\end{Proof}

\begin{Corollary}
\label{Corollaryfsplitparabolic}
Let $\kar\, k$ be a good prime for $G$. Then
the cotangent bundle $T^*(G/P)$ of $G/P$ is Frobenius split.
\end{Corollary}
\noindent
\begin{Proof}
This follows from Theorem \ref{Theoremparabolic} and Proposition
\ref{PropositionSpringer}.
\end{Proof}

\begin{Theorem}
Assume that $\kar\, k$ is a good prime for $G$. 
Let $\lambda\in X(P)$ be a $P$-regular weight. Then 
$$
\H^i(T^*(G/P), \pi^*\L(\lambda)) = \H^i(G/P, S\mathfrak{u}_P^*\otimes\lambda) =
0
$$
for $i>0$.
\end{Theorem}
\noindent
\begin{Proof} 
The proof follows \S \ref{Proofoftheorem}. One applies the Koszul resolution
for the short exact sequence of $P$-modules 
$0\rightarrow (\mathfrak{g}/\mathfrak{u}_P)^*\rightarrow
\mathfrak{g}^*\rightarrow \mathfrak{u}_P^*\rightarrow 0$. We get for $n\geq 1$ an
exact sequence
$$
\dots\rightarrow S^{n-1}\g^*\otimes \wedge^1(\g/\u_P)^*\otimes\lambda\rightarrow 
S^n \g^*\otimes \lambda\rightarrow S^n \u_P^*\otimes \lambda \rightarrow 0
$$
after tensoring with $\lambda$. Again the vanishing $\H^i(G/P, S\u_P^*\otimes\lambda)=0$
for any fixed $i>0$ follows from the vanishing
$$
\H^{i+j}(G/P, \wedge^j(\g/\u_P)^*\otimes\lambda) = 0
$$
for all $j\geq 0$. Since $\lambda$ induces an ample line bundle on $G/P$ this
vanishing follows when $\lambda$ is replaced by $n \lambda$ for all sufficiently large
$n$. In particular, we get the vanishing of $\H^i(T^*(G/P), \pi^*\L(p^r\lambda)) =
\H^i(G/P, S\u_P^*\otimes p^r \lambda)$ for any $i>0$ and all sufficiently large
$r$. 
Now the result follows from
Corollary \ref{Corollaryfsplitparabolic} and Lemma \ref{Lemmasplit}.
\end{Proof}

\section{The subregular nilpotent variety}

Throughout this section we assume that $G$ is simple (and simply connected) and
that $\kar\, k$ is good for $G$.

Let $\mathcal{U}$ be the unipotent variety in $G$ i.e. the closed subvariety of $G$
consisting of all unipotent elements. Then the map
$$
\varphi: G\times^B U \rightarrow \mathcal{U}
$$
mapping $(g, u)$ to $g u g^{-1}$ is a resolution of singularities 
(\cite{Humphreys}, Theorem 6.3) for all prime characteristics. If $P=P_{\alpha}$ is
a minimal parabolic subgroup associated with a short simple root $\alpha$, 
then $\varphi$ restricted to $G\times^B U_P$ factors through
$$
\varphi_\alpha: G\times^P U_P \rightarrow \mathcal{U}.
$$
\begin{Lemma}
\label{LemmaTits}
The map 
$$
\varphi_\alpha: G\times^P U_P \rightarrow S
$$
is birational onto its image $S$, which consists of
the variety of irregular elements (called the subregular unipotent variety).
\end{Lemma}
\noindent
\begin{Proof}
It follows by an argument of Tits that $\varphi_\alpha$ has connected fibres
(see \cite{Broer}, Proposition 4.2), so we need to show that $\varphi_\alpha$ is
separable. By Richardson's theorem (\cite{SpringerSteinberg}, I 5.1-5.6) 
the orbit maps for the conjugation action of $G$ on itself are separable for very 
good primes. This implies the separability of $\varphi_\alpha$ for good primes, when 
$G$ is not of type $A$. In type $A$ the separability of $\varphi_\alpha$ follows 
from the $\GL_n$-case, where the orbit maps for the conjugation action are 
separable for all primes.
\end{Proof}

By (\cite{BardsleyRichardson}, Corollary 9.3.4)  there is a (Springer) 
$G$-isomorphism between the 
unipotent variety $\mathcal{U}$ and the nilpotent cone $\mathcal{N}$ i.e. the closed subvariety of
$\g$ consisting of all nilpotent elements. In particular, we get that 
$\mathcal{N}$ is normal by the normality of $\mathcal{U}$ (\cite{Humphreys}, Proposition 1.3). Like in the 
unipotent case, the Springer resolution
$$
G\times^B \mathfrak{u} \rightarrow \mathcal{N}
$$
is a resolution of singularities, which gives a resolution (Lemma \ref{LemmaTits})
$$
\tilde{\varphi_\alpha}:G\times^P \mathfrak{n} \rightarrow \mathcal{S}
$$
of singularities of the subregular nilpotent variety $\mathcal{S}$, where $\mathfrak{n}$ is the
nilpotent radical of the Lie algebra of $P$. Let $\pi:T^*(G/B)\rightarrow G/B$ denote
the projection. 

\begin{Theorem}
The subregular nilpotent variety $\mathcal{S}$ is a normal Gorenstein variety with
rational singularies.
\end{Theorem}
\noindent
\begin{Proof}
The characteristic zero proof (\cite{Broer}, Theorem 4.4) carries over:
The closed subvariety $G\times^B \mathfrak{n}$ of the
cotangent bundle $G\times^B\mathfrak{u}$ is the zero scheme of 
a section of the pull back $\pi^*\L(-\alpha)$. So 
we get an exact sequence
$$
0\rightarrow \pi^*\L(\alpha)\rightarrow \O_{G\times^B\mathfrak{u}}
\rightarrow \O_{G\times^B \mathfrak{n}} \rightarrow 0.
$$
By Theorem \ref{Theoremvanishing} and the normality of $\mathcal{N}$, we get a 
short exact sequence
$$
0\rightarrow \H^0(T^*(G/B), \pi^*\L(\alpha))\rightarrow k[\mathcal{N}]
\rightarrow k[G\times^P\mathfrak{n}] \rightarrow 0.
$$
Let $\tilde{\mathcal{S}}$ denote the normalization of $\mathcal{S}$.
The surjection $k[\mathcal{N}]\rightarrow k[G\times^P \mathfrak{n}]$ factors through
the injection $k[S]\rightarrow k[\tilde{\mathcal{S}}]$ (followed by the 
map $k[\tilde{\mathcal{S}}]\rightarrow k[G\times^P \mathfrak{n}]$ induced by the 
normalization) 
via the restriction map
$k[\mathcal{N}]\rightarrow k[\mathcal{S}]$. This proves that $k[\mathcal{S}]=
k[\tilde{\mathcal{S}}]$ so that $\mathcal{S}$ is normal. 
By Theorem \ref{Theoremvanishing} the cohomology of $\O_{G\times^B \mathfrak{u}}$ vanishes. 
It follows that $\H^i(G\times^B \mathfrak{n}, \O_{G\times^B \mathfrak{n}}) = 
\H^i(G\times^P\mathfrak{n}, \O_{G\times^P \mathfrak{n}}) = 0$
for $i>0$, giving that  $\mathcal{S}$ has rational singularities (since $\tilde{\varphi_\alpha}$ 
is birational by Lemma \ref{LemmaTits}). As the canonical
line bundle of $G\times^P\mathfrak{n}$ is trivial, $\mathcal{S}$
is Gorenstein.
\end{Proof}

\section{Good filtrations}

Let $X$ be a smooth $B$-variety. A splitting section (or Frobenius splitting)
$\sigma\in\H^0(X, \omega_X^{1-p})$ is called canonical (\cite{Mathieu}, \cite{Kallen}, Definition 4.3.5) 
if $\sigma$ is $T$-invariant and for all $\alpha\in S$
$$
x_\alpha(t).\sigma = \sum_{i=0}^{p-1} t^i \sigma_{i, \alpha}
$$
for suitable $\sigma_{i, \alpha}\in \H^0(X, \omega_X^{1-p})$ (of weight $i\,\alpha$), where
$x_\alpha:k \rightarrow B$ is the root homomorphism corresponding to $\alpha$. 
 
Recall that a filtration 
$0=V_0\subset V_1\subset \dots$ of a $G$-module $V$ is called a good filtration if 
$V$ is the union of the $G$-submodules $V_0, V_1, \dots$
and $V_i/V_{i-1}\cong \H^0(G/B, \lambda_i)$ for $\lambda_i$ dominant.
We have the following weaker version of (\cite{Kallen}, Lemma 4.4.2) (due to Mathieu)

\begin{Lemma}
\label{LemmaKallen}
Let $X$ be a smooth $B$-variety and $\L$ a $G$-equivariant line bundle on 
$G\times^B X$. Assume that $G\times^B X$ admits a canonical splitting, then
the $G$-module $\H^0(G\times^B X, \L)$ has a good filtration.
\end{Lemma}

For good primes there is a $G$-equivariant map
$$
\varphi:\St\otimes\St \rightarrow \H^0(T^*(G/B), \O_{T^*(G/B)})
$$
such that $\varphi(a\otimes b)$ is a splitting section if
$\chi(a\otimes b) \neq 0$. Consider the splitting section of the
cotangent bundle $T^*(G/B)$ given by $\varphi(v^+\otimes v^-)$. 
It is easy to see that $\varphi(v^+\otimes v^-)$ is a canonical Frobenius splitting of 
$T^*(G/B)= G\times^B \mathfrak{u}$, since the definition can be checked for $v^+\otimes v^-\in
\St\otimes\St$.
\begin{Theorem}
\label{Theoremgoodfilt}
Suppose that $\kar\, k$ is a good prime for $G$. 
Let $\lambda\in X(T)$ be a weight (not necessarily dominant). Then
$$
\H^0(G/B, S^n \mathfrak{u}^*\otimes \lambda)
$$
has a good filtration for $n\geq 0$.
\end{Theorem}
\noindent
\begin{Proof}
By the above $T^*(G/B)=G\times^B \mathfrak{u}$ has a canonical Frobenius
splitting. This means that
$$
\H^0(T^*(G/B), \pi^*\L(\lambda)) = \H^0(G/B, S \mathfrak{u}^*\otimes \lambda)
$$
has a good filtration by Lemma \ref{LemmaKallen}, where $\pi: T^*(G/B)\rightarrow
G/B$ denotes the projection.
\end{Proof}

\begin{Remark}
Using Theorem \ref{Theoremparabolic} it follows in the same way that
$T^*(G/P) = G\times^P \u_P$ has a canonical Frobenius splitting for
a parabolic subgroup $P\supset B$. 
Mathieu has informed us that $H^0(X, \L)$ admits a good filtration 
if $X$ is a smooth $G$-variety with a canonical Frobenius splitting and
$\L$ a $G$-equivariant line bundle on $X$. In our case one may prove 
directly that
$G\times^B (G \times^P \u_P) \cong G/B \times (G\times^P \u_P)$
has a canonical Frobenius splitting, so that Lemma \ref{LemmaKallen} 
implies that $\H^0(G/P, S\u_P^*\otimes \lambda)$ has a good filtration for
(arbitrary) weights $\lambda\in X(P)$.
\end{Remark}

\begin{Theorem}
Suppose that $p>h$ and let $\lambda$ be a dominant weight. Then we have
an isomorphism for any $w\in W$ such that $w\cdot 0+ p\lambda$ is dominant
$$
\H^i(G_1, \H^0(G/B, w\cdot 0 + p\, \lambda))^{[-1]}\cong
\begin{cases}
\H^0(G/B, S^{(i-\ell(w))/2} \mathfrak{u}^*\otimes \lambda)& \text{if $i\equiv \ell(w) \mod 2$},\\
0 & \text{otherwise}.
\end{cases}
$$
where $()^{[-1]}$ denotes Frobenius (un)twist of a representation. 
In particular the cohomology of induced representations 
$\H^i(G_1, \H^0(G/B, w\cdot 0+ p\, \lambda))^{[-1]}$ admits a good
filtration.
\end{Theorem}
\noindent
\begin{Proof}
The key ingredient in the proof (in \cite{AndersenJantzen}) of the isomorphism is the
vanishing theorem \ref{Theoremvanishing}, which makes the spectral sequence 
(\cite{AndersenJantzen}, 3.3(2)) degenerate.
The good filtrations follow from
Theorem \ref{Theoremgoodfilt}.
\end{Proof}

\begin{Remark}
Andersen and Jantzen proved the above theorem for groups not having any components
of types $E$ and $F$(\cite{AndersenJantzen}, \S 5). For arbitrary $G$ they proved
the above theorem under the assumption that $\lambda$ is strongly dominant (\cite{
AndersenJantzen}, Corollary 3.7(b)).
\end{Remark}
\begin{Remark}
It follows from the linkage principle that the only dominant $\mu$ with 
$$
\H^\bullet(G_1, \H^0(G/B, \mu))\neq 0
$$ 
are of the form $w\cdot 0 + p\, \lambda$
for some $\lambda$ dominant and $w\in W$.
\end{Remark}

\section{Homogeneous Frobenius splittings}

The functions $(k[G]\otimes k[\u])^B = (k[G]\otimes S\u^*)^B$ on the cotangent bundle 
$T^*(G/B)$ have a natural grading. Let $\pi_d:(k[G]\otimes S\u^*)^B\rightarrow
(k[G]\otimes S^d\u^*)^B$ be the projection on the $d$-th homogeneous factor. Let $N$
denote the dimension of $G/B$. Then a function $f$ Frobenius splits $T^*(G/B)$
implies that $\pi_{N(p-1)}(f)$ Frobenius splits $T^*(G/B)$. A homogeneous
splitting function (of degree $N(p-1)$) descends to give a Frobenius splitting
of the projectivization $\P(T^*(G/B))$ (lines in $T^*(G/B)$) of the cotangent
bundle. These splittings are in some sense better behaved than the splittings
coming directly from $\St\otimes \St$.

\subsection{The $A_n$-case}

In type $A_n$ ($G=\SL_{n+1}(k)$) we have the $B$-equivariant isomorphism 
$\sigma: X\mapsto I+X$ between the
upper triangular nilpotent matrices $\u$ and the upper triangular unipotent
matrices $U$. In this way we see that the element $v^+\otimes v^-$ in $\St\otimes \St$
maps to the (splitting) function $f$
$$
(g, X) \mapsto \<v^+, g(X+I)g^{-1} v^-\>
$$
on the cotangent bundle $T^*(G/B) = G\times^B \u$ via $\H^0(\varphi)$ and $\sigma$. The 
function $g\mapsto \<v^+, g v^-\>$ is
a highest weight vector in $\St$ and equals the $(p-1)$-st power of the highest weight
function $f_\rho: g\mapsto \<w^+, g w^-\>$, where $w^+$ and $w^-$ are highest and lowest
weight vectors in $\H^0(G/B, \rho)$. The function $f_\rho$ is a product of certain highest
weight functions $f_{\omega_1}, \dots, f_{\omega_n}$, where $\omega_i$ denotes the
$i$-th fundamental dominant weight. Let $A = (a_{ij})_{1\leq ij\leq n+1}$ be a matrix in 
$G$, then it is well known that
$$
f_{\omega_s}(A) = \det((a_{ij})_{1\leq i,j \leq s})
$$
for $1\leq s\leq n$. In this way the (magical) splitting function of Mehta and van der Kallen 
\cite{MehtaKallen} on $T^*(G/B)$ is exactly $\pi_{N(p-1)}(f)$, where $N=n(n+1)/2$.
One interesting aspect of the Mehta - van der Kallen splitting 
is that it compatibly splits all $G\times^B \mathfrak{u}_P$, for any
parabolic subgroup $P\supseteq B$. Finding a suitable
splitting in this context for the other groups would be very interesting.


\begin{thebibliography}{99}

\bibitem{AndersenJantzen} Andersen, H.~H. and Jantzen, J.~C. 
Cohomology of induced representations of algebraic groups,
\emph{Math.~Ann.}, {\bf 269} (1984), 487--525

\bibitem{BardsleyRichardson} Bardsley, P. and Richardson, R.
\'Etale slices for algebraic tranformation groups in characteristic $p$,
\emph{Proc.~London~Math.~Soc.}, {\bf 51} (1985), 295--317

\bibitem{Broer} Broer, B.
Line bundles on the cotangent bundle of the flag variety,
\emph{Invent.~math.}, {\bf 113} (1993), 1--20

\bibitem{Broer1} Broer, B.
A vanishing theorem for Dolbeault cohomology of homogeneous vector bundles,
\emph{J.~reine~angew.~Math.}, {\bf 493}, (1997), 153--169

\bibitem{Broer2} Broer, B.
Normality of some nilpotent varieties and cohomology of line bundles on the 
cotangent bundle of the flag variety, 
\emph{Lie theory and geometry. In honor of Bertram Kostant} (eds. J.-L.
Brylinski, et al.), Birkh\"auser, Boston (1994), 1--19

\bibitem{Demazure} Demazure, M.
A very simple proof of Bott's theorem,
\emph{Invent.~math.}, {\bf 33} (1976), 271--272

\bibitem{Hartshorne} Hartshorne, R.
\emph{Algebraic Geometry},
Springer GTM 52

\bibitem{Humphreys} Humphreys, J.
\emph{Conjugacy Classes in Semisimple Algebraic Groups},
American Mathematical Society, Providence, 1995

\bibitem{Jantzen} Jantzen, J.~C.
\emph{Representations of Algebraic Groups},
Academic Press, Orlando, 1987


\bibitem{Kallen} van der Kallen, W.
\emph{Frobenius splittings and $B$-modules},
Tata Institute of Fundamental Research and Springer Verlag, Bombay 1993


\bibitem{LauritzenThomsen} Lauritzen, N. and Thomsen, J. F.
Frobenius splitting and hyperplane sections of flag manifolds,
\emph{Invent.~math.}, {\bf 128} (1997), 437--442 


\bibitem{Mathieu} Mathieu, O.
Filtrations of $G$-modules,
\emph{Ann.~Sci.~\'Ecole~Norm.~Sup.}, {\bf 23} (1990), 625--644

\bibitem{MehtaKallen} Mehta, V. and van der Kallen, W.
A simultaneous Frobenius splitting for closures of conjugacy
classes of nilpotent matrices,
\emph{Compos.~Math}, {\bf 84} (1992), 211--221

\bibitem{MehtaRamanathan} Mehta, V. and Ramanathan, A.
Frobenius splitting and cohomology vanishing for Schubert varieties, 
\emph{Annals~of~Math.}, {\bf 122} (1985), 27--40

\bibitem{Springer} Springer, T.
The unipotent variety of a semisimple group,
\emph{Algebraic Geometry} (ed.~S.~Abhyankar) 373--391, Oxford University Press, London 1969

\bibitem{SpringerSteinberg} Springer, T. and Steinberg, R.
Conjugacy classes,
\emph{Seminar on algebraic groups and related finite groups} (ed.~A.~Borel), 
Lecture Notes in Mathematics 131, Springer Verlag, New York 1970

\end{thebibliography}
\end{document}